\documentclass[a4paper,10pt]{article}
\usepackage{epsfig, graphicx, amssymb}
\usepackage{parskip}

\begin{document}

\title{
   \vspace{0cm}
   \textbf{Precise finite-sample quantiles of the\\
    Jarque-Bera adjusted Lagrange multiplier test\\[0.5cm]}
   }
\author{
   \emph{Diethelm W\"urtz and Helmut G. Katzgraber}\\[0.5cm]
   Swiss Federal Institute of Technology \\
   Institute for Theoretical Physics \\
   ETH H\"onggerberg, CH-8093 Z\"urich\\[0.5cm]
   \footnotesize\texttt{wuertz@itp.phys.ethz.ch, katzgraber@phys.ethz.ch}\normalsize\\[3cm]
}

\date{}

\maketitle

\vspace{-1cm}

\begin{center}
    First version: December 2004\\
    This version: August 2005
\end{center}

\vspace{1.5cm}

\begin{center}
\textbf{Abstract:}
\end{center}

\vspace{1cm} It is well known that the finite-sample null
distribution of the Jarque-Bera Lagrange Multiplier \emph{(LM)} test
for normality and its adjusted version \emph{(ALM)} introduced by
Urzua differ considerably from their asymptotic $\chi^2(2)$ limit.
Here, we present results from Monte Carlo simulations using $10^7$
replications which yield very precise numbers for the \emph{LM} and
\emph{ALM} statistic over a wide range of critical values and sample
sizes. This enables a precise implementation of the Jarque-Bera
\emph{LM} and \emph{ALM} test for finite samples.


\newpage

\section{Introduction}

The Jarque-Bera (1980, 1987) Lagrange multiplier test is likely
the most widely used procedure for testing normality of economic
time series returns. The algorithm provides a joint test of the
null hypothesis of normality in that the sample skewness $b_1$
equals zero and the sample kurtosis $b_2$ equals three. The null
is rejected when the Lagrange multiplier statistic

\vspace{0cm}

\begin{equation}
    LM = N \left (
        \frac{(b_1^{1/2})^2}{6}  +
        \frac{(b_2-3)^2}{24} \right )
\end{equation}

\vspace{0.1cm}

exceeds some critical value, which is taken in the asymptotic
limit from the $\chi^2(2)$ distribution. $N$ is the sample size,
$b_1^{1/2} = m_3/m_2^{3/2}$, $b_2 = m_4/m_2^2$ where $m_i$ is the
$i$-th central moment of the observations $m_i = \Sigma (x_j -
\overline{x})^i / N$, and $\overline{x}$ the sample mean.

\vspace{0.1cm}

Urzua (1996) modified the Jarque-Bera test replacing the
asymptotic means and variances by their exact finite-sample values
yielding

\begin{equation}
    ALM = N \left (
        \frac{(b_1^{1/2})^2}{c_1}  +
        \frac{(b_2-c_2)^2}{c_3} \right )
    ~.
\end{equation}

\vspace{0.1cm}

Here the parameters $c_{i}$, $i = 1$, $2$, $3$, are given by the expectation 
value ${\rm E}(\;)$ and variances ${\rm var}(\;)$ of the skewness and kurtosis

$$
    c_1 = {\rm var}(b_1^{1/2}) = \frac{6(N-2)}{(N+1)(N+3)} ~,
$$
$$
    c_2 = {\rm E}(b_2) = \frac{3(N-1)}{(N+1)} ~,
$$
$$
    c_3 = {\rm var}(b_2) = \frac{24N(N-2)(N-3)}{(N+1)^2(N+3)(N+5)} ~.
$$

\vspace{0.1cm}

Note, that the $ALM$ has the same asymptotic distribution as the
$LM$ statistic.

The work of Urzua (1996) as well as the work by Deb and Sefton
(1996) already warn about the incorrect use of the Jarque-Bera test
in the case of small- and medium-sized samples. The authors
have performed Monte Carlo simulations and have tabulated significance points
for 5\% and 10\%, on a series of sample sizes ranging between 10 and
800. Deb and Sefton used 600'000 replications in their Monte Carlo
simulations and Urzua used 10'000 replications and added results for
the 1\%, 15\% and 20\% significance points. Very recently Lawford
(2004) developed an accurate response surface approximation for the
5\% and 10\% critical values of the Jarque-Bera test based on Monte
Carlo simulations using 1 Million replications. The tables for the
$LM$ and $ALM$ statistics values presented in these papers are
restricted usually to a small set of parameters and the precision is
in most cases limited to two digits.

In this Letter we present tables with very precise values for both,
the $LM$ and $ALM$ statistic. Since the slow convergence of the
Monte Carlo simulation is well known we extend the simulations to 10
Million replications and enhance the mesh of $p$ values and sample
sizes considerably\footnote{The Monte Carlo data (about 1 Gbyte) are
available on request on 2 CDs from the authors.}.

The results have been used to implement R functions for the finite
sample Jarque-Bera test and the distribution itself, using either
the $LM$ or $ALM$ statistic. R (2004) is a powerful and widely used
GPL-licensed statistical software environment based on the S
language. In this sense our functions can also be called from the
commercial S-Plus software package. The R functions are part of the
Rmetrics software project, \emph{www.rmetrics.org}. The software is
GPL licensed, and the R package named \texttt{fSeries} can be
downloaded from the CRAN Server \emph{www.r-project.org}.


\section{Monte Carlo Simulation}

We have performed Monte Carlo simulations of the $LM$ and $ALM$
statistic using $10^7$ replications. The results are summarized in
Table 1 for both the $LM$ and $ALM$ statistic.

\begin{figure}[position=h]
  \begin{center}
  \vspace{0.2cm}
  \hspace{-0.2cm}
  \includegraphics[width=13.3cm]{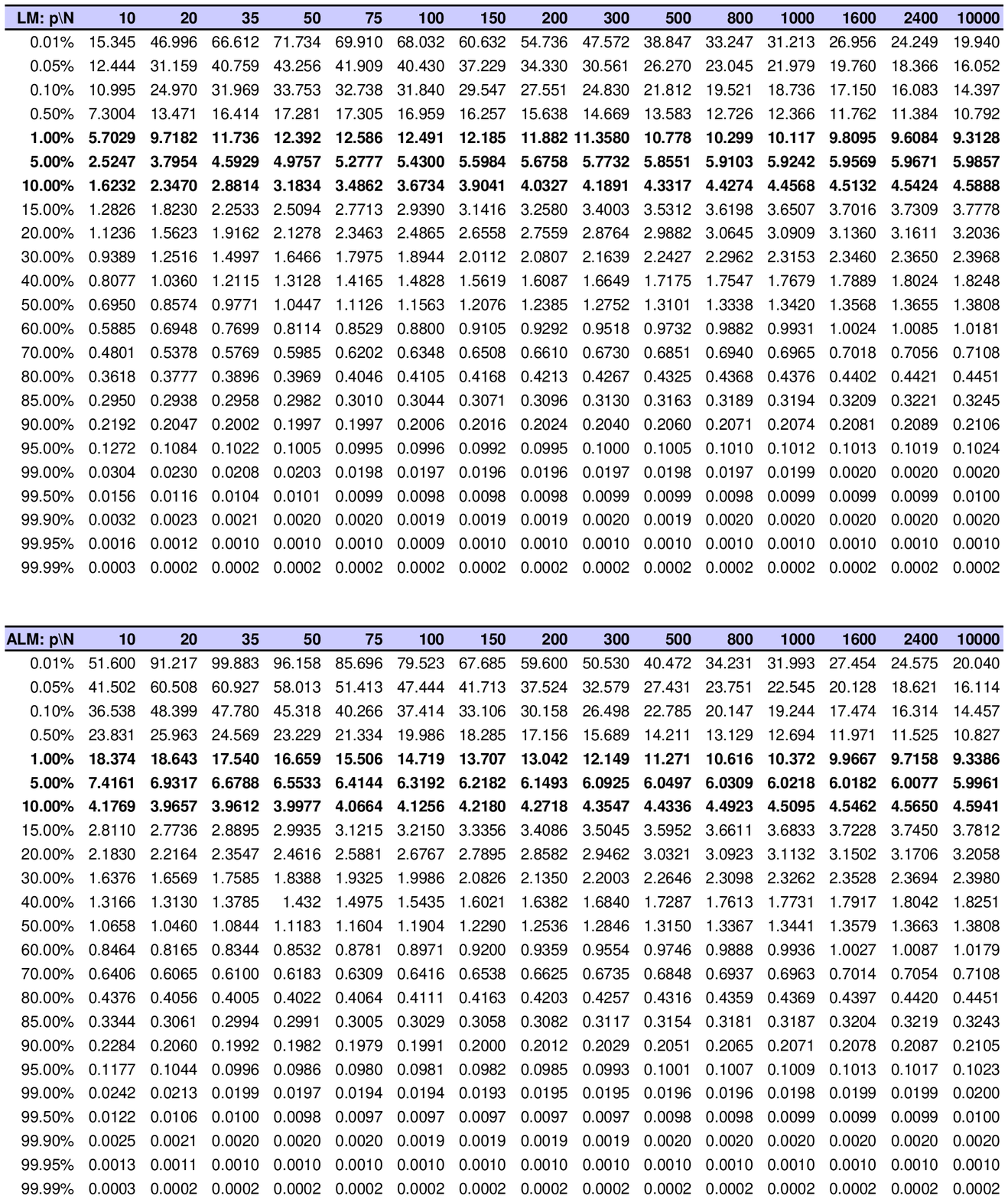}
  \vspace{0.3cm}
  \end{center}
  \footnotesize{\label{fig.jbTable}\emph{Table 1:
  Top: Significance points for the finite-sample Jarque-Bera
  test. Bottom: Same values for the adjusted Jarque-Bera Test.
  The numbers are based on Monte Carlo simulations using $10^7$
  replications. Note, that the $p$ values are listed in reverse
  order as $1-p$. The three major levels, 1\%, 5\% and 10\%,
  are written in  bold face.}}
\end{figure}

\vspace{0.2cm}

Figure 1 illustrates the results in a graph. The simulated $p$
values and the deviations from the asymptotic $\chi^2(2)$ limit are
shown. The curves belong to the same values of sample sizes $N$ as
listed in Table 1.

\newpage

\begin{figure}[t]
  \begin{center}
  \vspace{0.25cm}
  \hspace{-0.2cm}
  \includegraphics[width=15cm]{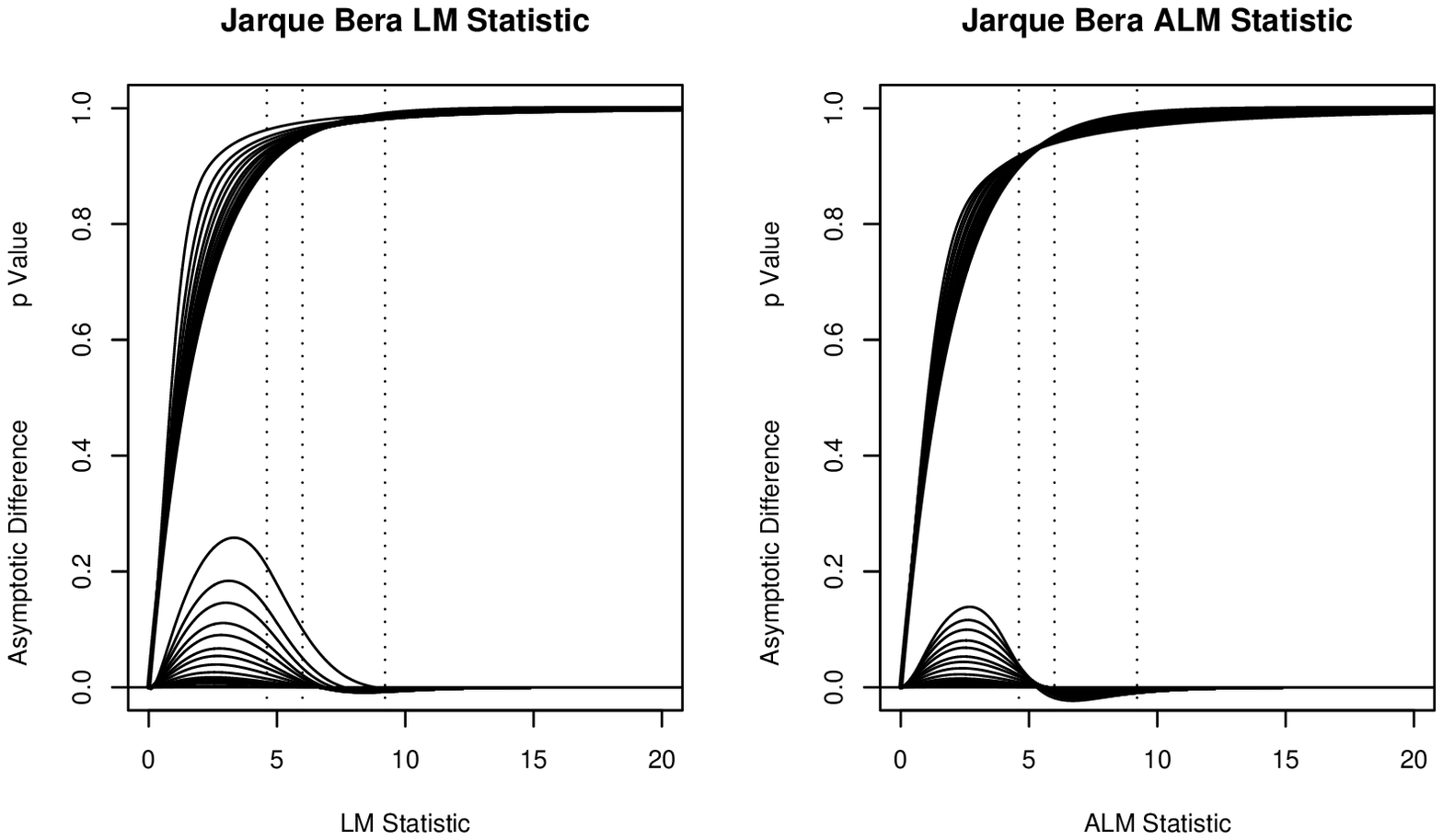}
  \vspace{0.2cm}
  \end{center}
  \footnotesize{\label{fig.jbFigure}\emph{Figure 1:
  $LM$ (left) and $ALM$ (right) finite-sample $p$ values and their
  differences with respect to the asymptotic limit. The upper bundle
  of curves shows the $p$ values. The lower bundle of curves
  measures the difference $p_N-p_\infty$ to the asymptotic limit.
  The graph clearly demonstrates that the adjusted Jarqua-Bera test
  outperforms the original version of the test. The three dotted
  vertical lines mark the 1\% (99\%), 5\% (95\%) and 10\% (95\%)
  levels in the asymptotic limit, respectively.}}
\end{figure}


\section{Response Surface}

To compute the $LM$ and $ALM$ statistic for a wide range of
quantiles and sample sizes one usually approximates the response
surface for a fixed value of $p$ as a series in powers of $1/N$

\vspace{-0.2cm}
\begin{equation}
    q(p, N) = q(p, \infty) + \sum_{k=1}^K \beta_k N^{-k} ~.
\end{equation}
\vspace{-0.2cm}

Lawford (2004) has done this for the 5\% and 10\% quantile lines. He
fitted his Monte Carlo data based on 1 Million replications for $K =
9$. The regression coefficients $\beta$ are listed in the
aforementioned paper. We have done fits over a wide range of
p-values. The results are shown in Figure 2 in comparison with those
obtained by Lawford. Note that Lawford's fit becomes less reliable
for small lengths where the convergence of the Monte Carlo
simulation slows down.

\begin{figure}[position=h]
  \begin{center}
  \vspace{0.3cm}
  \hspace{-0.3cm}
  \includegraphics[width=15cm]{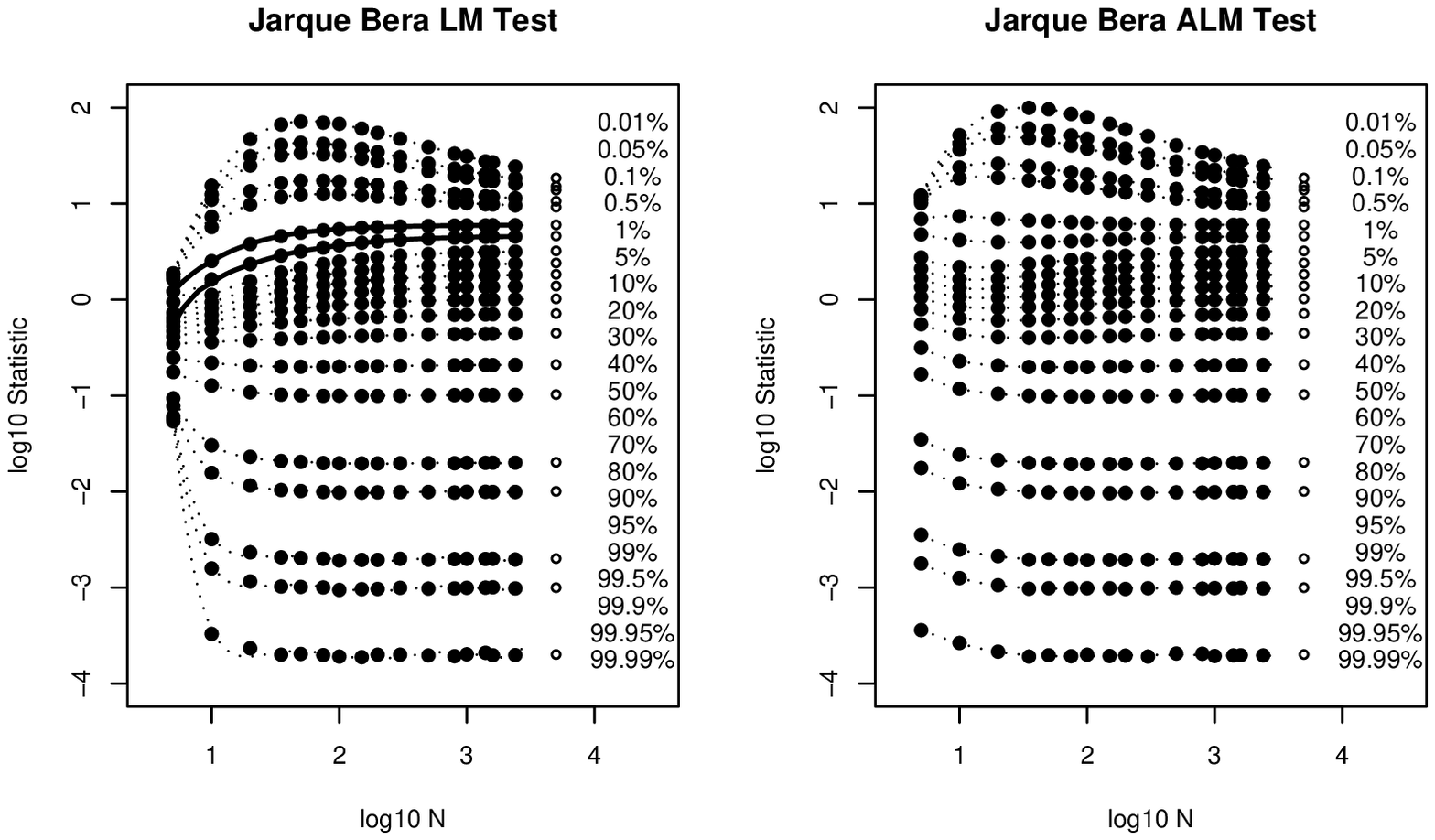}
  \vspace{0.15cm}
  \end{center}
  \normalsize{\label{fig.jbQuantiles}\emph{Figure 2:
  LM (left) and ALM (right) statistic for
  a wide range of $p$ values as a function of sample sizes. The
  dots show the results from the Monte Carlo simulations
  using $10^7$ replications together with the asymptotic limit
  (marked by the open circles). The dotted lines are fitted
  series expansions of order $K=6$ in $1/N$. The two thick
  lines in the left LM graph display the results of Lawford
  for the 5\% and 10\% levels.}}
\end{figure}

Another approach would be an Edgeworth (1917) expansion of the
distribution in $1/N$. Unfortunately, we have found that the
expansion converges extremely slowly. Thus we have applied ``Curve Fitting'',
as suggested by Rothenberg (1984), to approximate the response
surface. Simple linear interpolation, 2-dimensional splines or
connectionist function approximators are only three possibilities
from many others. We have followed the first approach fitting on
logarithmic scales. The results are shown in Figure 3 for both the
traditional Jarque-Bera test as well as its adjusted version.

\begin{figure}[position=h]
  \begin{center}
  \vspace{0.3cm}
  \hspace{-0.2cm}
  \includegraphics[width=14cm]{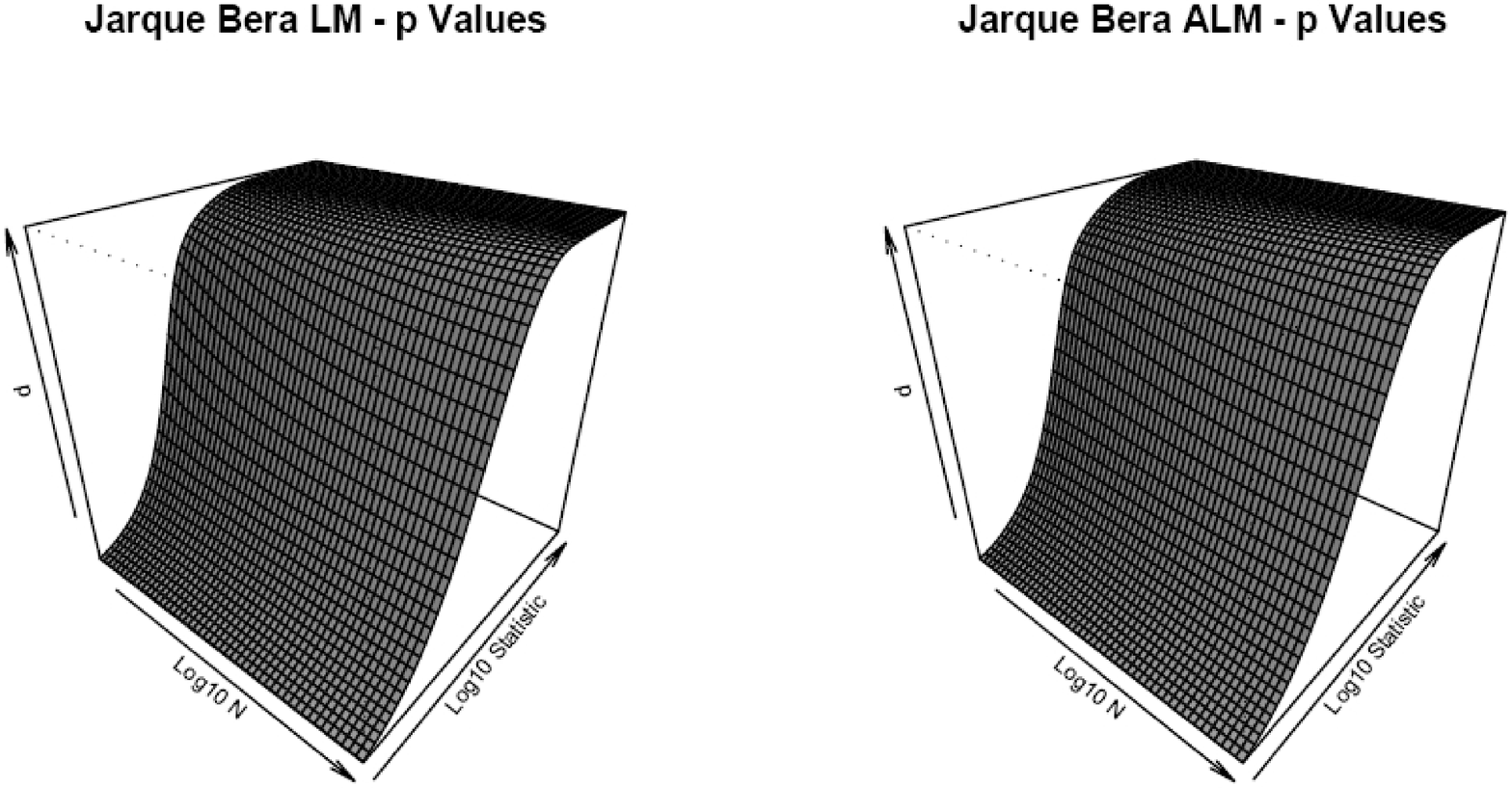}
  \vspace{0.3cm}
  \end{center}
  \normalsize{\label{fig.jbSurface}\emph{Figure 3:
  LM (left) and ALM (right) surface
  of $p$ values for a wide range of statistics (0.4 ... 100) and
  sample sizes (10 ... 10'000). Note, that the x- and y-axis are
  on logarithmic scales. The inputs consist of almost 2000
  p-values ranging between $0.0001$ and $0.9999$.}}
\end{figure}

We have implemented the Jarque-Bera test for finite samples into
an R package named \texttt{JarqueBera}. The underlying simulations
with $10^7$ replications have been done with a separate C program using
a multiplicative lagged Fibonacci random number generator with a
lag of size 1279. The software comes with the following functions:

\vspace{-0.2cm}
\begin{verbatim}
            pjb(q, N = Inf, method = c("LM", "ALM"))
            qjb(p, N = Inf, method = c("LM", "ALM"))

            jb.test(x)
\end{verbatim}

The first two functions compute the distribution function
\texttt{pjb()} and the quantile function \texttt{qjb()} either for
the $LM$ (default option) or $ALM$ test version. By default the
asymptotic values from the $\chi^2(2)$ distribution with two degrees
of freedom are returned. If \texttt{N} is specified by an integer,
the finite-sample values will be returned. The arguments \texttt{x},
and \texttt{q} are vectors of quantiles, and \texttt{p} is a vector
of probabilities. \texttt{n} denotes the number of observations. The
function \texttt{jbTest()} computes for a sample of empirical
observations given by the vector \texttt{x} the $LM$ and the $ALM$
statistic and returns the associated $p$ values for the finite
sample and in the asymptotic limit.

In what follows, we illustrate the use of the above R
functions. We use a simulated series drawn from a Student-t
distribution, \texttt{rt()}, which becomes normal in the limit of an infinite
number of degrees of freedom and becomes more and more fat tailed
when the number of degrees of freedom is decreased. Typical values
as can be found from economic time series are in the range of
about 3 to 4.

\vspace{0.5cm}
\begin{verbatim}
            > jb.test(x = rt(n = 100, df = Inf))

                    Jarque-Bera Test

            data:  rt(n = 100, df = Inf)
            LM = 1.9333, ALM = 1.883,
            LM p-value = 0.291, ALM p-value = 0.323,
            p-value = 0.3804


            > jb.test(x = rt(100, df = 4))

                    Jarque-Bera Test

            data:  rt(100, df = 4)
            LM = 1933.967, ALM = 2239.254,
            LM p-value = NA, ALM p-value = NA,
            p-value < 2.2e-16
\end{verbatim}
\vspace{0.5cm}

The Jarque-Bera test function \texttt{jbTest()} returns the $LM$ and
$ALM$ statistic and computes the associated finite-sample
$p$ values. In addition, the asymptotic $p$ value is printed for
comparison. If the finite-sample $p$ values are too small
\texttt{NA} will be returned. This is the case for the fat tailed
Student-t distribution in the above example where the p-value is of
the order of the $10^{-16}$.


\vspace{0.3cm}

\section{Summary}

This Letter tabulates precise $p$ values for the Jarque-Bera finite
sample normality test. In addition to the original version of the
Lagrange Multiplier test we have also computed finite-sample
$p$ values for its adjusted version formulated by Urzua (1996). In
contrast to previous investigations the results are derived from a
MC simulation with $10^7$ replications. To our knowledge this is one
of the largest simulations ever done in statistics. The outcome of
the simulation are very precise values for finite samples which we
have tabulated and can now be used for an improved hypothesis
testing. The test can be used together with the statistical R or
S-Plus environments. It should be straightforward to implement the
test procedure also into other statistical software packages like
Matlab, Eviews, SAS, or others.


\vspace{0.3cm}

\section*{References}


\hspace{-0.5cm}

\noindent
Bowman K., Shenton L.,
    \emph{Omnibus test contours for departures from normality
    based on $b_1$ and $b_2$},
    Biometrika 62, 1975, 243--250.

\noindent
Deb P., Sefton M.,
    \emph{The distribution of a Lagrange multiplier test of
    normality},
    Economics Letters 51, 1996, 123--130.

\noindent
Edgeworth F.Y.,
    \emph{On the mathematical representation of statistical data},
    Journal of the Royal Statistical Society 80, , 1917, 411--437.

\noindent
Jarque C.M, Bera A.K.,
    \emph{Efficient tests for normality, homoscedasticity and
    serial independence of regression residuals},
    Economics Letters 6, 1980, 255--259.

\noindent
Jarque C.M, Bera A.K.,
    \emph{A test for normality of observations and regression
    residuals},
    International Statistical Review 55, 1987, 163--172.

\noindent
Lawford S.,
    \emph{Finite-sample quantiles of the Jarque-Bera test},
    Brunel University Preprint, 2004, 7 pages.

\noindent
R Core Team,
    \emph{R Manuals},
    downloadable from: http://cran.r-project.org.

\noindent Rmetrics,
    \emph{Teaching Financial Engineering and Computational Finance with R},
    http://www.rmetrics.org.

\noindent
Rothenberg T.J.,
    \emph{Approximating the distributions of econometric
    estimators and test statistics},
    Handbook of Econometrics, Volume II, 1984, 881--935.

\noindent
Urzua M.,
    \emph{On the correct use of omnibus tests for normality},
    Economics Letters 53, 1996, 247--251.

\end{document}